\documentclass[12pt]{amsart}
\usepackage{amssymb,latexsym,amsmath,amscd,amsthm,graphicx, color}
\usepackage[all]{xy}
\usepackage{pgf,tikz}
\usepackage{mathrsfs}
\usepackage{breqn}
\usepackage{fancyhdr}
\usepackage[numbers]{natbib}
\usepackage{setspace}
\usepackage[english]{babel}
\usepackage[autostyle]{csquotes}
\usepackage{verbatim}
\DeclareFontFamily{U}{mathx}{\hyphenchar\font45}
\DeclareFontShape{U}{mathx}{m}{n}{
      <5> <6> <7> <8> <9> <10>
      <10.95> <12> <14.4> <17.28> <20.74> <24.88>
      mathx10
      }{}
\DeclareSymbolFont{mathx}{U}{mathx}{m}{n}
\DeclareFontSubstitution{U}{mathx}{m}{n}
\DeclareMathAccent{\widecheck}{0}{mathx}{"71}
\DeclareMathAccent{\wideparen}{0}{mathx}{"75}

\usetikzlibrary{arrows}
\definecolor{uuuuuu}{rgb}{0.26666666666666666,0.26666666666666666,0.26666666666666666}
\definecolor{xdxdff}{rgb}{0.49019607843137253,0.49019607843137253,1.}
\definecolor{ffqqqq}{rgb}{1.,0.,0.}

\newtheorem{theorem}{Theorem}[section]
\newtheorem{lemma}[subsection]{Lemma}

\theoremstyle{definition}
\newtheorem{remark}[subsection]{Remark}
\newtheorem{definition}[subsection]{Definition}

\theoremstyle{remark}

\newtheorem{proposition}[subsection]{Proposition}
\numberwithin{equation}{section}
\def\ni{\noindent}

\begin{document}
\pagenumbering{gobble}

\title[]{A note on boundedness of singular integral operators on functions spaces over Local Fields}

\author{Salman Ashraf and Qaiser Jahan}






\keywords{Local fields, Besov spaces, Triebel–Lizorkin spaces, Hardy spaces, Singular integral operators.}

\begin{abstract}
In this article we consider the classical singular integral operator over a local field with rough kernels. We study the boundedness of such an operator on different function spaces by relaxing the smoothness condition on kernels.
\end{abstract}

\maketitle

\section{\textbf{Introduction}} 

In the fundamental paper \cite{CA}, Singular integral theory were introduced by A. P. Calder\'on and A. Zygmund, which plays an important role in Harmonic analysis. Let $\Omega : \mathbb{R}^n \rightarrow \mathbb{R}$ be a Lebesgue measurable function with $\int_{\mathbb{S}^{n-1}} \Omega(y)dy = 0.$ We say that a linear operator $T_\Omega : S(\mathbb{R}^n) \rightarrow S^\prime(\mathbb{R}^n)$ is a rough singular integral operator if 
\begin{align*}
T_\Omega f(x) = \lim_{\epsilon \to 0} \int_{|y| > \epsilon} \dfrac{\Omega(y/|y|)}{|y|^n}f(x-y) dy.
\end{align*} 

The rough singular integral operators are one of the main topics in the theory of singular integral operators. In a series of remarks below we will see some known results about $L^p$ boundedness of $T_\Omega.$

\begin{itemize} \setlength\itemsep{.5em}
\item[(\romannumeral 1)] In 1956, Calder\'on and Zygmund \cite{CA1} proved that if $\Omega \in L\log^+L(\mathbb{S}^{n-1})$ is even or $\Omega \in L^1(\mathbb{S}^{n-1})$ is odd then $T_\Omega$ is of type $(p,p)$ for $1 < p < \infty.$
\item[(\romannumeral 2)] In 1979, $L^p$ boundedness of $T_\Omega$ for $1 < p < \infty$ proved independently by Connett \cite{CW} and Ricci-Weiss \cite{RG} under less restrictive condition that $\Omega \in H^1(\mathbb{S}^{n-1})$ since 
\begin{align*}
L\log^+L(\mathbb{S}^{n-1}) \subset H^1(\mathbb{S}^{n-1}).
\end{align*} 
Here $H^1(\mathbb{S}^{n-1})$ is the Hardy space on the unit sphere.
\item[(\romannumeral 3)] In 1998, Grafakos and Stefanov in \cite{GS}, considered the condition 
\begin{align} \label{1}
\sup_{\xi \in \mathbb{S}^{n-1}} \int_{\mathbb{S}^{n-1}} \big|\Omega(\theta)\big|\bigg(\log\dfrac{1}{|\theta\cdot\xi|}\bigg)^{1+\alpha}d\theta ~< ~ +\infty 
\end{align}
for $\alpha > 0$ and gave $L^p$ boundedness of $T_\Omega$ if $\Omega$ satisfies (\ref{1}). Also in \cite{GS1}, Grafakos and Stefanov gave a survey of known results in the theory of Calder\'on-Zygmund singular integral operator $T_\Omega.$ 
\end{itemize}

The boundedness of $T_\Omega$ on homogeneous Besov space $\dot{B}^s_{p,q}$ and on homogeneous Triebel-Lizorkin space $\dot{F}^s_{p,q}$ were studied by Chen, Fan and Ying in \cite{CFY}. The result of \cite{CFY} is as follows:

\begin{theorem} (See \cite{CFY}) Let $1 < p,q < \infty$ and $s \in \mathbb{R}.$ If $\Omega \in L^r(\mathbb{S}^{n-1}),~r > 1$ with mean value zero, then 
\begin{itemize}
\item[(\romannumeral 1)] $\|T_\Omega f\|_{\dot{F}_{p,q}^s(\mathbb{R}^n)}  \leq C \|f\|_{\dot{F}_{p,q}^s(\mathbb{R}^n)}$ 
\item[(\romannumeral 2)] $\|T_\Omega f\|_{\dot{B}_{p,q}^s(\mathbb{R}^n)}  \leq C  \|f\|_{\dot{B}_{p,q}^s(\mathbb{R}^n)}.$ 
\end{itemize}
\end{theorem}

The result of \cite{CFY} then improved by Chen and Zhang in \cite{CZ}:

\begin{theorem} (See \cite{CZ}) Let $1 < p,q < \infty$ and $s \in \mathbb{R}.$ If $\Omega \in L^1(\mathbb{S}^{n-1})$ with mean value zero and satisfies
\begin{align*}
\sup_{\xi \in \mathbb{S}^{n-1}} \int_{\mathbb{S}^{n-1}} \big|\Omega(\theta)\big|\bigg(\log\dfrac{1}{|\theta\cdot\xi|}\bigg)^{1+\alpha}d\theta ~< ~ +\infty \\
\end{align*}
for all $\alpha > 1$ then $\|T_\Omega f\|_{\dot{F}_{p,q}^s(\mathbb{R}^n)}  \leq C \|f\|_{\dot{F}_{p,q}^s(\mathbb{R}^n)}.$ 
\end{theorem}
 
Chen and Ding in \cite{CD}, considered the condition that $\Omega \in H^1(\mathbb{S}^{n-1})$ and proved the following result:  

\begin{theorem} (See \cite{CD}) Let $1 < p,q < \infty$ and $s \in \mathbb{R}.$ If $\Omega \in H^1(\mathbb{S}^{n-1})$ with mean value zero, then 
\begin{itemize}
\item[(\romannumeral 1)] $\|T_\Omega f\|_{\dot{F}_{p,q}^s(\mathbb{R}^n)}  \leq C \|f\|_{\dot{F}_{p,q}^s(\mathbb{R}^n)}$ 
\item[(\romannumeral 2)] $\|T_\Omega f\|_{\dot{B}_{p,q}^s(\mathbb{R}^n)}  \leq C  \|f\|_{\dot{B}_{p,q}^s(\mathbb{R}^n)}.$ \\
\end{itemize}
\end{theorem}

In this paper, our main aim is to establish the boundedness of singular integral operator with Hardy space kernel on  Lebesgue spaces, Triebel–Lizorkin spaces, and Besov spaces over local fields. In \cite{MZ}, the boundedness of $p$-adic singular integral operators and their commutators were obtained on $p$-adic Morrey spaces. In \cite{QY}, $\delta$-type Calderón–Zygmund operators were introduced and established the boundedness on weighted Hardy spaces on locally compact Vilenkin groups. Recently, KP Ho \cite{HK} extended the studies of $\delta$-type Calderón–Zygmund operators on Morrey and Hardy-Morrey spaces in locally compact Vilenkin groups.

A local field $K$ is a locally compact, totally disconnected, non-Archimedian norm valued and non-discrete topological field, we refer \cite{MHT} to basic Fourier analysis on local fields. The local fields are essentially of two types (excluding the connected local fields $\mathbb{R}$ and $\mathbb{C}$) namely characteristic zero and of positive characteristic. The characteristic zero local fields include the $p$-adic field $\mathbb{Q}_p$ and the examples of positive characteristic are the Cantor dyadic groups, Vilenkin $p$-groups and $p$-series fields.

Besov spaces $B^s_{rt}(K)$ and Triebel-Lizorkin spaces $F^s_{rt}(K)$ on local fields were introduced and studied by Onneweer and Su \cite{Os}. In \cite{book_impulse4}, Weiyi Su, introduced, \textquote{$p$-type smoothness} of the functions defined on local fields. To complete the theoretical base of the Function spaces on local fields and to broaden the range of its applications, a series of studies have been carried out such as the construction theory of functions on local fields \cite{book_impulse10}, the Weierstrass functions, Cantor functions and their $p$-adic derivatives in local fields, the Lipschitz classes on local fields etc. (see \cite{book_impulse5}). Wavelet theory on local fields has also been developed widely. Jiang, Li and Jin \cite{JLJ} have introduced the concept of multiresolution analysis and wavelet frames on local fields of positive characteristic \cite{JL}. Later, Behera and Jahan have developed the theory of wavelets on such a field (see \cite{BJ1,BJ4}). Recently, Lukomskii and Vodolazov \cite{LV} gave the fast algorithm for the discrete Fourier transform on local fields of characteristic zero.

The classical singular integral operators of Calderón–Zygmund type over local fields were considered by Phillips and  Taibleson (see \cite{PT,PK,MHT}). Let $K$ be a local field and $\mathfrak{D}^* = \{\xi \in K : |\xi|=1\}.$ Let $\Omega \in L^\infty(K)$  and satisfy

\begin{align} \label{2}
\Omega(\mathfrak{p}^jx) &= \Omega(x)~ \text{for all}~j \in \mathbb{Z}, \\ \label{3}
\int\limits_{\mathfrak{D^{*}}} \Omega(x)dx &= 0 .
\end{align}

The singular integral operator is defined by 

\begin{align} \label{4}
Tf(x) = \lim\limits_{k \to -\infty} \int_{|y| > q^{k}} f(x-y) \dfrac{\Omega(y)}{|y|} dy,~~~~ \text{for}~ k \in \mathbb{Z},
\end{align}

and the truncated singular integral operator is defined as

\begin{align} \label{5}
T_kf(x) =  \int_{|y| > q^{k}} f(x-y) \dfrac{\Omega(y)}{|y|} dy,~~~~ \text{for}~ k \in \mathbb{Z},
\end{align}

initially for $ f \in S(K).$ In particular, Taibleson \cite{MHT} proved the following result:

\begin{theorem} \label{t1}
Suppose that $\Omega \in L^{\infty}(K)$ satisfies the conditions (\ref{2})-(\ref{3}) and 

\begin{align} \label{6}
\sup_{|y|=1} \sum\limits_{j=1}^{\infty} \int_{\mathfrak{D}^*}~ |\Omega(x+\mathfrak{p}^jy) - \Omega(x)|~ dx < \infty.
\end{align}

Then $T_k$ is bounded in $L^r(K)$ for $1 < r < \infty.$ Furthermore $T(f) = \lim_{k\to-\infty} T_k(f)$ exists in the $L^r(K)$ norm.
\end{theorem}

In the series of remarks above, we have seen that boundedness of $T_\Omega$ is still possible on $L^p(\mathbb{R}^n),$ homogeneous Besov space $(\dot{B}^s_{p,q})$ and homogenous Triebel-Lizorkin spaces $(\dot{F}^s_{p,q})$ without any smoothness on $\Omega.$ In the setting of Local field $K$, as in Theorem \ref{t1}, Taibleson considered the smoothness condition (\ref{6}) on $\Omega,$  which is analogues of H\"ormander's condition and gave the boundedness of $T_k$ on $L^r(K)(1 < r < \infty)$. Our interest is in defining the truncated singular integral operator $T_k$ with rough kernel and therefore we replaces the smoothness condition (\ref{6}) with $\Omega \in H^1(\mathfrak{D}^*)$ and studied the boundedness of $T_k$ on the Lebesgue spaces, Triebel–Lizorkin spaces, and Besov spaces over Local fields. In particular, we proved the following results:

\begin{theorem} \label{thm1}
Let $ 1 < r < \infty,~ k \in \mathbb{Z} $ and $T_k$ be the operator defined by (\ref{5}). If $\Omega \in H^1(\mathfrak{D}^*)$ satisfies (\ref{2}) and (\ref{3}), then 
\begin{align*}
\|T_kf\|_{L^r(K)} \leq Cq^{-k}\|\Omega\|_{H^1(\mathfrak{D}^*)} \|f\|_{L^r(K)},
\end{align*}

where $C > 0$ is independent of $k.$ 
\end{theorem}

\begin{theorem} \label{thm2}
Let $1 < r, t < \infty,~s > 0,~ k \in \mathbb{Z} $ and $T_k$ be the operator defined by (\ref{5}). If $\Omega \in H^1(\mathfrak{D}^*)$ satisfies (\ref{2}) and (\ref{3}), then  

\begin{itemize}
\setlength\itemsep{.5em}
\item[(\romannumeral 1)] $\|T_kf\|_{F_{r t}^s(K)}  \leq Cq^{-k} \|\Omega\|_{H^1(\mathfrak{D}^*)} \|f\|_{F_{r t}^s(K)},$ 
\item[(\romannumeral 2)] $\|T_kf\|_{B_{r t}^s(K)}  \leq Cq^{-k} \|\Omega\|_{H^1(\mathfrak{D}^*)} \|f\|_{B_{r t}^s(K)},$ 
\end{itemize}

where $C > 0$ is independent of $k.$
\end{theorem}

We remark here that by Theorem \ref{thm1} and Theorem \ref{thm2}, it is possible to have boundedness of $T_k$ with $\Omega \in H^1(\mathfrak{D}^*)$ on $L^r(K), B_{r t}^s(K)$ and $F_{r t}^s(K)$ but the operator norm is dependent on $k$ by constant $q^{-k}$ and therefore convergence of $T_k$ as $k \rightarrow - \infty$ is not possible in the respective function spaces norm, which is not the case when  $\Omega$ satisfies smoothness condition (\ref{6}) as in Theorem \ref{t1}.


This paper is organized as follows. In Section 2, we provide a brief introduction to local fields and test function class and distribution on such a field. We have also provided the definition of $p$-type derivatives, Triebel-Lizorkin spaces $F^s_{rt}(K)$ and Besov spaces $B^s_{rt}(K)$ on local fields. Further, we have given the definition of atomic Hardy space on $\mathfrak{D}^*$ and some Lemmas for the proof of Theorem \ref{thm1} and Theorem \ref{thm2}. Section 3 and section 4 are devoted to the proof of Theorem \ref{thm1} and Theorem \ref{thm2}, respectively.  

\section{Preliminary definitions and lemmas}

\subsection{Local Fields}
Let $K$ be a field and a topological space. If the additive group $K^{+}$ and multiplicative group $K^{*}$ of $K$ both are locally compact Abelian group, then $K$ is called \textit{locally compact topological field}, or \textit{local field.}

If $K$ is any field and is endowed with the discrete topology, then $K$ is a local field. Further, if $K$ is connected, then $K$ is either $\mathbb{R}$ or $\mathbb{C}$. If $K$ is not connected, then it is totally disconnected. So, here local fields means a field $K$ which is locally compact, non-discrete and totally disconnected. If $K$ is of characteristic zero, then K is either a $p$-adic field for some prime number $p$ or a finite algebraic extension of such a field. If $K$ is of positive characteristic, then $K$ is either a field of formal Laurent series over a finite field of characteristic $p$ or an algebraic extension of such a field.

Let $K$ be a local field. Since $K^{+}$ is a locally compact abelian group, we choose a Haar measure $dx$ for $K^{+}$. If $ \alpha \neq 0,~ \alpha \in K$, then $d(\alpha x)$ is also a Haar measure and by the uniqueness of Haar measure we have $d(\alpha x)= |\alpha|dx$. We call $|\alpha|$ the \textit{absolute value} or \textit{valuation} of $\alpha$ which is, in fact, a natural non-Archimedian norm on $K$. A mapping $|\cdot|~ : K \rightarrow \mathbb{R}$ satisfies:

\begin{itemize}
\item[(a)] $|x|=0$ if and only if $x=0$;
\item[(b)] $|xy|= |x||y|$ for all $x,y \in K$;
\item[(c)] $|x+y| \leq \text{max} \{|x|, |y|\},$  for all $x,y \in K$, \
\end{itemize}

The property $(c)$ is called the \textit{ultrametric inequality}. It follows that
\begin{equation*}
|x+y| = \text{max} \{|x|, |y|\}~ \text{if}~ |x| \neq |y|.
\end{equation*}

The set $\mathfrak{D}= \{x \in K : |x| \leq 1\}$ is called the ring of integers in $K$. It is the unique maximal compact subring of $K$. The set $\mathfrak{P}=\{x \in K : |x| < 1\}$ is called the prime ideal in $K$. It is the unique maximal ideal in $\mathfrak{D}$.  It can be proved that $\mathfrak{D}$ is compact and open and hence $\mathfrak{P}$ is also compact and open. Therefore, the residue space $\mathfrak{D}/\mathfrak{P}$ is isomorphic to a finite field $GF(q)$, where $q =p^c$ for some prime number $p$ and $c \in \mathbb{N}$.

Let $\mathfrak{p}$ be a fixed element of maximum absolute value in, $\mathfrak{P}$ and it is called a \textit{prime element} of $K$. For a measurable subset $E$ of $K$, let $|E|=\int_K\bold{1}_E(x)dx$, where $\bold{1}$ is the indicator function of $E$ and $dx$ is the Haar measure of $K$ normalized so that $|\mathfrak{D}|=1$. Then it is easy to prove that $|\mathfrak{P}| = q^{-1}$ and $|\mathfrak{p}| = q^{-1}$ (see~\cite{MHT}).

Let $\mathcal{U}=\{a_i\}_{i=0}^{q-1}$ be any fixed full set of coset representative of $\mathfrak{P}$ in, $\mathfrak{D}$ then each $x \in K$ can be written uniquely as
$x =\sum\limits_{l=k}^{\infty}c_l\mathfrak{p}^k$, where $c_l\in\mathcal{U}$.
If $x (\neq 0) \in K$ then $|x| =q^k$ for some $ k \in \mathbb{Z}.$

Let $\mathfrak{D}^{*} = \mathfrak{D}\setminus \mathfrak{P} = \{x \in \mathfrak{D} : |x|=1\}$; $\mathfrak{D}^*$ is the group of units in $K^*$. If $x \neq 0$, we can write as $x = \mathfrak{p}^k x^\prime$, with $ x^{\prime}\in  \mathfrak{D}^*$. Let $\mathfrak{P}^k = \mathfrak{p}^k \mathfrak{D} = \{x \in K : |x| \leq q^{-k}\}~, k \in \mathbb{Z}$. These are called \textit{fractional ideals}.

The set $\{\mathfrak{P}^k \subset K: k \in \mathbb{Z}\}$ satisfies the following:

\begin{itemize}
\item[(\romannumeral 1)] $\{\mathfrak{P}^k \subset K : k \in \mathbb{Z}\}$ is a base for neighborhood system of identity in $K$, and $\mathfrak{P}^{k+1} \subset \mathfrak{P}^k,~ k \in \mathbb{Z};$
\item[(\romannumeral 2)] $\mathfrak{P}^k ,~ k \in \mathbb{Z},$ is open, closed and compact in $K$;
\item[(\romannumeral 3)] $K = \bigcup\limits_{{k=-\infty}}^{+\infty} \mathfrak{P}^k~~\text{and}~~\{0\}=\bigcap\limits_{{k=-\infty}}^{+\infty} \mathfrak{P}^k.$
\end{itemize} 

Let $\Gamma$ be the character group of the additive group $K^{+}$ of $K$. There is a nontrivial character $\chi \in \Gamma$ which is trivial on $\mathfrak{D}$ but is non-trivial on $\mathfrak{P}^{-1}$. Let $\chi$ be a non-trivial character on $K^{+}$, then the corresponding relationship $\lambda \in K \longleftrightarrow \chi_{\lambda} \in \Gamma$ is determined by $\chi_{\lambda}(x)=\chi(\lambda x),$ and the topological isomorphism is established for $K$ and $\Gamma$, moreover, we have $\Gamma=\{ \chi_{\lambda} : \lambda \in K \}.$

For $k \in \mathbb{Z}$, let $\Gamma^{k}$ be the annihilator of $\mathfrak{P}^{k},$ that is,

\begin{align*}
 \Gamma^{k} = \{\chi \in \Gamma : \forall x \in \mathfrak{P}^{k} \implies \chi(x)=1\}.
\end{align*}

The locally compact group $\Gamma=\{ \chi_{\lambda} : \lambda \in K \}$ can be endowed with the non-Archimedian norm as $ \Gamma^{k}=\{ \chi_{\lambda} \in \Gamma : |\lambda| \leq q^{k}\}=\mathfrak{p}^{-k}\mathfrak{D},~k \in \mathbb{Z}.$

For the Haar measure $dx$ of $K^{+}$, let $d\xi$ be the Haar measure on $\Gamma$, chosen such that 

$$|\Gamma^{0}|=|\mathfrak{D}|=1, ~\text{and}~|\Gamma^{k}| = \dfrac{1}{|\mathfrak{P}^{k}|} = q^{k}.$$

The set $\{\Gamma^k \subset \Gamma :~ k \in \mathbb{Z}\}$ satisfies the following properties:

\begin{itemize}
\item[(\romannumeral 1)] $\Gamma^{k} \subset \Gamma^{k+1},~ k \in \mathbb{Z}$, are increasing sequence in $\Gamma$ which is open, closed and compact;
\item[(\romannumeral 2)] $\Gamma = \bigcup\limits_{{k=-\infty}}^{+\infty} \Gamma^k~~\text{and}~~\{I\}=\bigcap\limits_{{k=-\infty}}^{+\infty} \Gamma^k,$ where $I$ is the unit element of $\Gamma$;
\item[(\romannumeral 3)] $\{\Gamma^k\}_{k\in \mathbb{Z}}$ is a base of the unit $I$ of character group $\Gamma$.
\end{itemize}

\subsection{Function classes}

In this subsection we will introduce some function classes:

\textit{Test function class $S(K)$:} It is the linear space in which functions have the form

$$ \phi(x)= \sum_{j=1}^{n} c_j  \Phi_{\mathfrak{P}^{j}}(x-h_j),~c_j \in \mathbb{C}, h_j \in K, j \in \mathbb{Z},n \in \mathbb{N},$$ 

where $\Phi_{\mathfrak{P}^{j}}$ is the characteristic function of $\mathfrak{P}^{j}.$ The space $S(K)$ is an algebra of continuous functions with compact support that separates points. Consequently, $S(K)$ is dense in $C_0(K)$ as well as $L^r(K),~1 \leq r < \infty.$ Similarly, the test function class $S(\Gamma)$ on $\Gamma$ can be defined. However, since $K$ is isomorphic to $\Gamma$, so $S(K)$ and $S(\Gamma)$ can be regarded as equivalent with respect to absolute value or valuation.

The space $S(K)$ is equipped with a topology as a topological vector space as follows: Define a null sequence in $S(K)$ as a sequence $\{\phi_n\}$ of functions on $S(K)$ in such a way that each $\phi_n$ is constant on cosets of $\mathfrak{P}^l$ and is supported on $\mathfrak{P}^k$ for a fixed pair of integers $k$ and $l$ and the sequence converges to zero uniformly. The space $S(K)$ is complete and separable and is called the \emph{space of testing function}.

Since $S(K)$ is dense in $L^1(K)$, thus the Fourier transformation of $\phi(x) \in S(K)$ is defined by

$$\widehat{\phi}(\xi) \equiv (\mathcal{F}\phi)(\xi) = \int_{K} \phi(x)\bar{\chi}_\xi(x)dx,~~~  \xi  \in \Gamma,$$

and the inverse Fourier transformation of $\phi \in S(K)$ is defined by the formula 

$$ \widecheck{\phi}(x) \equiv (\mathcal{F}^{-1}\phi)(x) = \int_{\Gamma} \phi(\xi)\chi_x(\xi)d\xi,~~~  x  \in K.$$

$S^{\prime}(K)$, the \textit{space of distributions}, is a collection of continuous linear functional on $S(K)$. $S^{\prime}(K)$ is also a complete topological linear space. The action of $f$ in $S^{\prime}(K)$ on an element $\phi$ in $S(K)$ is denoted by $\langle f, \phi \rangle$. The distribution space $S^{\prime}(K)$ is given the weak$\ast$ topology. Convergence in $S^{\prime}(K)$ is defined in the following way: $f_k$ converges to $f$ in $S^{\prime}(K)$ if  $\langle f_k, \phi \rangle$ converges to $\langle f,\phi\rangle$ for any $\phi \in S(K).$

The Fourier transformation $\widehat{f}$ of a distribution  $f \in S^{\prime}(K)$ is defined by

$$ \langle \widehat{f}, \phi \rangle ~ = ~ \langle f, \widehat{\phi} \rangle~~~\text{for~all}~ \phi \in S(K),$$ 

 the inverse Fourier transformation $\widecheck{g}$ is defined by 

$$ \langle \widecheck{g}, \psi \rangle = \langle g, \widecheck{\psi} \rangle~~~\text{for~all}~ \psi \in S(K).$$

\subsection{Derivatives and Integrals} 
In \cite{book_impulse4}, Weiyi Su introduced $p$-type derivatives and $p$-type integrals for general locally
compact Vilenkin group using pseudo-differential operators, which also includes derivatives and integrals of fractional orders.

\begin{definition} 
Let $\alpha > 0$, if for a complex Haar measurable function $f : K \rightarrow \mathbb{C}$, the integral

$$T_{\langle\cdot\rangle^\alpha} f(x) \equiv \int_{\Gamma} \bigg\{\int_{K} \langle \xi \rangle^\alpha f(t) \bar{\chi}_\xi(t-x)dt \bigg\} d\xi$$

\ni exists at $x \in K$, then $T_{\langle\cdot\rangle^\alpha} f(x)$ is said to be a point-wise $\alpha$-order $p$-type derivative of $f(x)$ at $x$, and it is denoted by $f^{\langle\alpha\rangle}(x).$
\end{definition}

\begin{definition} Let $\alpha > 0$, if for a complex Haar measurable function $f : K \rightarrow \mathbb{C}$, the integral

$$T_{\langle\cdot\rangle^{-\alpha}} f(x) \equiv \int_{\Gamma} \bigg\{\int_{K} \langle \xi \rangle^{-\alpha} f(t) \bar{\chi}_\xi(t-x)dt \bigg\} d\xi$$

\ni exists at $x \in K$, then $T_{\langle\cdot\rangle^{-\alpha}} f(x)$ is said to be a point-wise $\alpha$-order $p$-type integral of $f(x)$ at $x$, denoted by $f_{\langle\alpha\rangle}(x).$
\end{definition}

\begin{remark} We suppose that $ \alpha > 0$ in the above definitions, thus the order $\alpha$ of $p$-type derivatives and integrals can be any positive real numbers, thus, fractional order derivatives and fractional order integrals all are contained. Moreover, for $\alpha =0$, 
$$f^{\langle 0 \rangle}(x) = f(x) = f_{\langle 0 \rangle}(x)~~\text{for}~x\in K.$$
\end{remark}

\subsection{Function spaces on Local fields}

\begin{definition} \label{def2} 
 Let $S(K)$ and $S(\Gamma)$ be the test function classes on $K$ and on $\Gamma$, respectively, $S^{\prime}(K)$ and $S^{\prime}(\Gamma)$ are their distribution spaces, respectively.

Take a sequence $\{\phi_j(\xi)\}_{j=0}^{+\infty} \subset S(\Gamma)$ satisfying 

\begin{itemize}
\item[(\romannumeral 1)]   
$\text{supp}~\phi_0 \subset \{\xi \in \Gamma  : |\xi|< q\}$,\\
$\text{supp}~\phi_j \subset \{\xi \in \Gamma : q^{j-1} < |\xi|< q^{j+1}\},~~~~j \in \mathbb{N}$
 
\item[(\romannumeral 2)] 
$\sum_{j=0}^{\infty} \phi_j(\xi)=1,~~\xi \in \Gamma.$

\item[(\romannumeral 3 )] $\bigg|(\widecheck{\phi}_j)^{\langle s \rangle}(x)\bigg| \leq c_s q^{-j+js},~s \in (0, +\infty),~~~j \in \mathbb{N}_0,~~x \in K,$

\ni where $(\widecheck{\phi}_j)^{\langle s \rangle}$ is the  point-wise $s$-order $p$-type derivatives of $\widecheck{\phi}_j$ and $c_s$ is a constant depending on $s$ only.
\end{itemize}

We denote  
$$ A(\Gamma)= \bigg\{ \{\phi_j\}_{j\in \mathbb{N}_0} \subset S(\Gamma)~\text{satisfying}~ (\romannumeral 1) - (\romannumeral 3)\bigg\}.$$
\end{definition}

Specially, for a local field $K$, we just need to take $\{\phi_j\}_{j=0}^{\infty} \subset A(\Gamma)$ as follows:
\begin{align*}
\phi_0(\xi) & = \Phi_{\Gamma^0}(\xi),\\
\phi_j(\xi) & = \Phi_{\Gamma^j \setminus \Gamma^{j-1}}(\xi),~~~j \in \mathbb{N},
\end{align*}
where $\Phi_{A}$ is the characteristic function of $A$. It satisfies the following conditions
\begin{itemize}
\item[(\romannumeral 1)] $\text{supp}~\Phi_{\Gamma^0} \subset \{\xi \in \Gamma  : |\xi|<q\}$,\\
$\text{supp}~\Phi_{\Gamma^j \setminus \Gamma^{j-1}} \subset \{\xi \in \Gamma : q^{j-1} < |\xi|< q^{j+1}\},~~~~j \in \mathbb{N},$\\

\item[(\romannumeral 2)] $ \Phi_{\Gamma^0}(\xi) + \sum_{j=1}^{\infty} \Phi_{\Gamma^{j}\setminus \Gamma^{j-1}}(\xi)=1,~~~~\xi \in \Gamma,$\\

\item[(\romannumeral 3 )] $\bigg|(\widecheck{\Phi}_{\Gamma^0})^{\langle s \rangle}(x)\bigg| \leq c_s,$\\
$\bigg|(\widecheck{\Phi}_{\Gamma^{j}\setminus \Gamma^{j-1}})^{\langle s \rangle}(x)\bigg| \leq c_s q^{-j+js},~s \in (0, +\infty),~~~j \in \mathbb{N},~~x \in K,$\\~\text{where}~$c_s$~\text{is a constant depending on} $s$ \text{only}.
\end{itemize}

\begin{definition} \label{def1}
Let $\{ \phi_j\}_{j=0}^{\infty} = \{\Phi_{\Gamma^0},~\Phi_{\Gamma^{j}\setminus \Gamma^{j-1}}\}_{j=1}^{\infty} \subset A(\Gamma),$ we define the operators $\Delta_j$ as follows
\begin{align*}
\Delta_j f =  \mathcal{F}^{-1} \Big(\phi_j \mathcal{F}f\Big),~~~~~~~~~~j \in \mathbb{N}_0,~~f \in S^{\prime}(K).
\end{align*}

Then we obtain the Littlewood-Paley decomposition,

$$f = \sum_{j=0}^{\infty} \Delta_j f,$$

of all $f \in S^{\prime}(K).$
\end{definition}

Note that $\Delta_j f \in S^{\prime}(K)$ for any $f \in S^{\prime}(K)$.

\begin{proposition} 
Let $u$ be in $S^{\prime}(K).$ Then, we have, $u = \sum_{j=0}^{\infty} \Delta_j u,$ in the sense of the convergence in the space $S^{\prime}(K).$
\end{proposition}

\begin{definition} \textbf{(B-type space)}  
For $0 < r \leq +\infty,~0 < t \leq +\infty,~ s \in \mathbb{R},$ we define B-type spaces or Besov spaces on local fields $K$ as

$$ B_{rt}^s(K)=\{f \in S^{\prime}(K):~ \|f\|_{{ B_{rt}^s}(K)} < \infty\}, $$

\ni with norm

\begin{align*}
\|f\|_{{ B_{rt}^s}(K)} &= \|q^{sj} \Delta_j f \| _{\ell_t(L^r(K))}\\
&= \Bigg\{\sum\limits_{j=0}^{\infty} \|q^{sj} \Delta_j f\|^t_{L^r(K)}   \Bigg\}^{\frac{1}{t}}\\
&= \Bigg\{\sum_{j=0}^{\infty} q^{sjt} \bigg\{ \int_{K} |\Delta_j f|^r dx \bigg\}^{\frac{t}{r}}  \Bigg\}^{\frac{1}{t}}.
\end{align*}
\end{definition}

\begin{definition} \textbf{(F-type space)} 
For $0 < r < \infty,~0 < t \leq \infty,~ s \in \mathbb{R},$ we define F-type spaces or Triebel-Lizorkin spaces  on local field $K$ as
$$ F_{rt}^s(K)=\{f \in S^{\prime}(K):~ \|f\|_{{ F_{rt}^s}(K)} < \infty\}, $$

\ni  with norm

\begin{align*}
\|f\|_{{ F_{rt}^s}(K)} &= \|q^{sj} \Delta_j f \| _{L^r(\ell^t(K))}\\
&= \Bigg\|\Bigg\{\sum_{j=0}^{\infty} |q^{sj} \Delta_j f|^t   \Bigg\}^{\frac{1}{t}} \Bigg\|_{L^r(K)}\\
&= \Bigg( \int_{K} \Bigg\{\sum_{j=0}^{\infty}q^{sjt} | \Delta_j f|^t   \Bigg\}^{\frac{r}{t}} dx \Bigg)^{\frac{1}{r}}.
\end{align*} 
\end{definition}

\begin{remark} The spaces $A_{rt}^s(K),~ A \in \{B,F\}$ are independent of the selection of sequences  $\{\phi_j\}_{j\in \mathbb{N}_0} \subset A(\Gamma).$ They are quasi-Banach spaces (Banach spaces for $r,t \geq 1$), and $ S(K) \subset A_{rt}^s(K) \subset S^{\prime}(K),$ where the first embedding is dense if $r < \infty$ and $t < \infty.$ The theory of the spaces $A_{rt}^s(K)$  has been developed  in \cite{book_impulse5}.
\end{remark}




Now, we give the definition of atomic Hardy space on $\mathfrak{D}^*,$ which can be obtained by the idea in \cite{DJE}.


\begin{definition} \label{def1}
($H^1(\mathfrak{D}^*)$ atoms) A function $a: \mathfrak{D}^{*} \rightarrow \mathbb{C}$ is
a $(1, \infty)$ atom if 

\begin{itemize}
\setlength\itemsep{.5em}
\item[(\romannumeral 1)] $\text{supp}~ a \subset \mathfrak{D}^{*} ,$
\item[(\romannumeral 2)] $ \|a\|_{L^\infty(K)}~ \leq~ |\mathfrak{D}^{*}|^{-1} = (1-q^{-1})^{-1},$
\item[(\romannumeral 3)] $\int_{\mathfrak{D}^{*}} a(x)dx ~=~0.$ 
\end{itemize}
\end{definition}

\begin{definition}
If $f \in H^1(\mathfrak{D}^*)$ then there exists a sequence $\{\lambda_i\}_{i \in \mathbb{N}} \in \ell^1$ and a sequence of $(1, \infty)$ atoms $\{a_i\}_{i \in \mathbb{N}}$ such that $f = \sum\limits_{i=1}^{\infty} \lambda_ia_i$ and $\|f\|_{H^1(\mathfrak{D}^*)} \approx \sum\limits_{i=1}^{\infty} |\lambda_i|.$
\end{definition}

\begin{theorem} \cite{BJ4} (Minkowski's integral inequality) Let $(X, \mu)$ and $(Y, \nu)$ be $\sigma-$finite measure spaces. Let $1 \leq r < \infty $ and $F$ be a function on the product space $X \times Y.$ Then 

\begin{align} \label{7}
\bigg(\int_X \bigg(\int_Y |F(x,y)| d\nu(y)\bigg)^r d\mu(x)\bigg)^{1/r} ~ \leq ~ \int_Y \bigg(\int_X |F(x,y)|^r d\mu(x)\bigg)^{1/r} d\nu(y). 
\end{align}
\end{theorem}

We will use this inequality for the counting measure $\mu = \sum_{k=0}^{\infty} \delta_k $ on $X = \mathbb{N}_0:$
\begin{align} \label{8}
\bigg(\sum_{k=0}^{\infty}\bigg(\int_Y |F_k(y)| d\nu(y)\bigg)^r \bigg)^{1/r} ~ \leq ~ \int_Y \bigg(\sum_{k=0}^{\infty} |F_k(y)|^r\bigg)^{1/r} d\nu(y).
\end{align}


The following Lemmas are useful for us and given in \cite{MHT}.
\begin{lemma} \label{lemma1} \setstretch{1.20}
 There are constants $A,~B > 0$ depending only on $K$ such that if $f \in L^1(K),~ f(x) \geq 0, ~ \lambda > 0$ then there is a countable collections of mutually disjoints spheres $\{W_i\}_{i \in \mathbb{N}}~ ,~D_\lambda = \cup_i W_i$, and a decomposition of $f = f_1 + f_2$ such that:
 
\begin{itemize} \setlength\itemsep{.5em}
\item[(\romannumeral 1)] $|D_\lambda|~ = ~\sum\limits_{i} |W_i| < \|f\|_{L^1} \lambda^{-1}$ 
\item[(\romannumeral 2)] $|f(x)|~ \leq~ \lambda~~~\text{a.e. for}~ x \notin D_\lambda$ 
\item[(\romannumeral 3)] $|f_2(x)| ~\leq ~q \lambda~~~\text{a.e. for}~ x \in D_\lambda$ 
\item[(\romannumeral 4)] $f_2(x)~ =~ f(x),~  x \notin D_\lambda$ 
\item[(\romannumeral 5)] $f_1(x)~ = ~0,~  x \notin D_\lambda$ 
\item[(\romannumeral 6)] $\int\limits_{W_i} f_1(x)dx = 0~~~\forall i.$
\end{itemize}
\end{lemma}

\begin{remark} \label{rem1} The conclusions of above lemma yield the following results:

\begin{itemize}
\setlength\itemsep{.5em}
\item[(\romannumeral 1)] $f_1 \in L^1(K),~~~\|f_1\|_{L^1(K)} ~\leq ~\|f\|_{L^1(K)}$
\item[(\romannumeral 2)] $f_2 \in L^1(K),~~~\|f_1\|_{L^1(K)} ~= ~\|f\|_{L^1(K)}$
\item[(\romannumeral 3)] $f_2 \in L^\infty(K),~~~\|f_2\|_{L^\infty(K)} ~\leq~  q\lambda$ 
\item[(\romannumeral 4)] $f_2 \in L^2(K),~~~\|f_2\|^2_{L^2(K)}~\leq ~ \lambda q \|f\|_{L^1}.$
\end{itemize}
\end{remark}

\begin{lemma} \label{lemma2} \setstretch{1.20}
Suppose $T$ is a linear operator defined on $S(K)$ with values in $S^{\prime}(K)$. Suppose further that when $f \in S(K)$ then $Tf \in L^2(K) $ and $\|Tf\|_{L^2(K)} \leq A_2 \|f\|_{L^2(K)},~A_2 > 0$ independent of $f$ and $Tf$ is a function such that for all $\lambda > 0,~~ |\{x \in K: |Tf(x)| >  \lambda \}|  \leq A_1 \|f\|_{L^1(K)} \lambda^{-1},A_1 > 0$ independent of $f$ and $\lambda.$ Then there is a constant $A_r > 0$ independent of $f,~ 1 < r \leq 2$ such that for $f \in S(K),$ $Tf \in L^r(K)$ and $\|Tf\|_{L^r(K)} \leq A_r \|f\|_{L^r(K)}.$ 
\end{lemma}

\begin{definition} An operator $T$ that is defined on $S(K)$ with values in $S^{\prime}(K)$ is said to be commute with translations if for all $h \in K,$ $\tau_h T = T_{\tau_h}.$ 
\end{definition}

\begin{lemma} \label{lemma3} \setstretch{1.20} 
Suppose $T$ is a linear operator defined on $S(K)$ with values in $S^{\prime}(K)$ that commutes with translations. Suppose further that $\phi \in S(K)$ implies that $T\phi \in L^r(K),~ \|T\phi\|_{L^r(K)} \leq A \|\phi\|_{L^r(K)}, ~ A > 0$ independent of $\phi,~ 1 \leq  r \leq \infty.$ Let $ \frac{1}{r} + \frac{1}{r^\prime} =1.$ Then $T\phi \in L^{r^\prime}(K)$ for all $\phi \in S(K)$ and $\|T\phi\|_{L^{r^\prime}(K)} \leq A \|\phi\|_{L^{r^\prime}(K)}.$   
\end{lemma}

\section{Proof of Theorem \ref{thm1}}

We begin by discretising $T_k.$ Let $\Omega \in H^1(\mathfrak{D}^*)$ and for $k \in \mathbb{Z}$ we have 

\begin{align} \nonumber
T_k(f)(x) & =  \int\limits_{|y| > q^{k}} f(x-y) \dfrac{\Omega(y)}{|y|} dy, \\ \nonumber
& =  \int\limits_{|y| = q^{k+1}} f(x-y)  \dfrac{\Omega(y)}{|y|} dy +  \int\limits_{|y| =q^{k+2}} f(x-y)  \dfrac{\Omega(y)}{|y|} dy + \cdots, \\ \nonumber
& = \sum\limits_{j = k }^{\infty }~ \int\limits_{|y| =q^{j+1}} f(x-y) \dfrac{\Omega(y)}{|y|}dy, \\ \label{eq1}
& =  \sum\limits_{j = k }^{\infty }q^{-(j+1)}~ \int\limits_{|y| =q^{j+1}} f(x-y) \Omega(y)dy. \\ \nonumber
\end{align}

Since $\Omega \in H^1(\mathfrak{D}^*)$ then we can express $\Omega$ as 

\begin{equation} \label{eq2}
\Omega = \sum_{i=1}^{\infty} \lambda_ia_i, 
\end{equation}

where $\{\lambda_i\}_{i \in \mathbb{N}} \in l^1$ and $a_i \textquotesingle s$ are $(1, \infty)$ atoms. Putting the value of (\ref{eq2}) into (\ref{eq1}), we get 

\begin{align} \nonumber
T_k(f)(x) & = \sum\limits_{i=1}^{\infty}\sum\limits_{j = k }^{\infty }q^{-(j+1)}\lambda_i \int\limits_{|y| =q^{j+1}} f(x-y) a_i(y)dy, \\ \label{eq3}
& = \sum\limits_{i=1}^{\infty} \lambda_i B_if,
\end{align}
where  
\begin{align*}
B_if = \sum\limits_{j = k }^{\infty }q^{-(j+1)}~ \int\limits_{|y| = q^{j+1}} f(x-y) a_i(y)dy. 
\end{align*}

For simplicity, we denote $B_i$ and $a_i$ by $B$ and $a$, respectively. Define 
\begin{align*}
g_j(y) = a(y) \Phi_{A_j}(y),
\end{align*}

where $\Phi_{A_j}$ is the characteristic function of $A_j = \{y~ : ~ |y|=q^{j+1}\}.$ Then 

\begin{align} \label{eq4}
Bf(x) = \sum\limits_{j = k }^{\infty }q^{-(j+1)} g_j \ast f(x).
\end{align}

From (\ref{eq3}) we get

\begin{align} \nonumber
\|T_kf\|_{L^r(K)} & \leq  \sum_{i=1}^{\infty} |\lambda_i| \|Bf\|_{L^r(K)} , \\ \label{eq5}
& \leq C^\prime \|\Omega\|_{H^1(\mathfrak{D}^*)} \|Bf\|_{L^r(K)}. 
\end{align}

In order to prove that $T_k$ is bounded on $L^r(K)~(\text{for}~1 < r < \infty)$, it is enough to show the same for $B.$ We will do this proof in two steps. In the first step we will show that  $B$ is bounded on $L^2(K).$ In second step we will prove that $B$ is weak type $(1, 1)$ and conclude from Lemma \ref{lemma2} that $B$ is bounded in $L^r(K)$ for $1 < r \leq 2.$ Then we will use Lemma \ref{lemma3} to obtain the boundedness of $B$ in $L^r(K)$ for $2 < r < \infty.$ 

\subsection{$L^2(K)$ boundedness of $B$} Let $f \in L^2(K)$ then $g_j \ast f \in L^2(K),$ since $g_j \in L^1(K).$ Also

\begin{align*} 
|\widehat{g_j}(\xi)| & = \bigg| \int_{K} a(x)\Phi_{A_j}(x)\bar{\chi}_\xi(x)dx \bigg|, \\
&= \bigg| \int_{\mathfrak{D}^*} a(x)\bar{\chi}_\xi(x)dx \bigg|,~~~~~~~~~~(\because~ \text{supp}~ a\Phi_{A_j} \subset \mathfrak{D}^*) \\ 
& \leq \int_{\mathfrak{D}^*} |a(x)||\bar{\chi}_\xi(x)|dx,\\ 
& \leq \|a\|_{L^\infty(K)} ~\cdot~ |\mathfrak{D}^*| \\
& \leq |\mathfrak{D}^{*}|^{-1} |\mathfrak{D}^*| = 1.  
\end{align*}

By Plancherel, we have that 
\begin{align*}
\|g_j \ast f\|_{L^2(K)} & \leq ~\|\widehat{g_j} \|_{L^{\infty}(K)}  \| \widehat{f}\|_{L^2(K)}\\
& \leq ~  \| \widehat{f}\|_{L^2(K)}, 
\end{align*}

\ni and

\begin{align} \nonumber
\|B f\|_{L^2(K)} & = \bigg\|\sum_{j = k }^{\infty } q^{-(j+1)}~  g_j \ast f \bigg\|_{L^2(K)}\\ \nonumber
& \leq ~ \sum_{j = k }^{\infty } q^{-(j+1)}~ \|g_j \ast f\|_{L^2(K)} \\ \nonumber
&\leq ~ \sum_{j = k }^{\infty } q^{-(j+1)}~ \| f\|_{L^2(K)} \\ \label{eq6}
&= ~  q^{-k}(q-1)^{-1} \| f\|_{L^2(K)} .
\end{align}




\subsection{Weak $(1, 1)$ boundedness of $B$}  Let $f \in L^1(K)$ and fix $\lambda > 0.$ We apply the decomposition in lemma \ref{lemma1} to obtain the functions $f_1,~f_2$ so that $f = f_1 + f_2.$ Now observe that 

\begin{align*}
\{x \in K :  |Bf(x)| >  \lambda\} \subset \{x \in K : |Bf_1(x)| >  \lambda/2\} \bigcup \{x \in K : |Bf_2(x)| >  \lambda/2\}  \equiv E_1 \cup E_2,
\end{align*}

\ni and so, we have 

\begin{align} \label{eq7}
|\{x \in K :  |Bf(x)| >  \lambda\}|~  \leq ~ |\{x \in K : |Bf_1(x)| >  \lambda/2\}| + |\{x \in K : |Bf_2(x)| >  \lambda/2\}|
\end{align}
and need to estimates each of these terms. The estimates of $E_2$ is easy since we have

\begin{align} \nonumber
|E_2| = \int_K \Phi_{E_2}(x)dx =\int_{E_2} 1~ dx & ~\leq ~\dfrac{4}{{\lambda}^2} \|B f_2(x)\|^2_{L^2(K)} \\ \nonumber
& ~\leq ~ \dfrac{4}{{\lambda}^2}   \| f_2\|^2_{L^2(K)} \\ \nonumber
& ~\leq~ \dfrac{4}{{\lambda}^2}   \lambda q\| f\|_{L^1(K)} \\ \label{eq8}
& ~= ~\dfrac{4q }{\lambda}\|f\|_{L^1(K)}. 
\end{align}

Here, we have used that $B$ is bounded on $L^2(K)$ and property (\romannumeral 4) in Remark \ref{rem1}. Now we find estimate on $E_1.$ We have

\begin{align} \nonumber
|E_1| ~& =~ |\{x \in K : |Bf_1(x)| >  \lambda/2\}|  \\ \label{eq9}
&= ~\bigg|\bigg(\bigcup_iW_i\bigg)\bigcap \{x \in K : |Bf_1(x)| >  \lambda/2\} \bigg|~ + ~ \bigg|\bigg(\bigcup_iW_i\bigg)^\complement   \bigcap \{x \in K : |Bf_1(x)| >  \lambda/2\} \bigg|. 
\end{align}

Consider the first term of (\ref{eq9}), 

\begin{align} \nonumber
\bigg|\bigg(\bigcup_iW_i\bigg)\bigcap \{x \in K : |Bf_1(x)| > \lambda/2\} \bigg| & ~=~ \bigg|\bigcup_iW_i\cap \{x \in K : |Bf_1(x)| >  \lambda/2\} \bigg| \\ \nonumber
& ~\leq~ \bigg|\bigcup_iW_i\bigg| \\ \nonumber
&~ \leq ~\sum\limits_i |W_i| \\ \label{eq10}
&~ \leq~ \|f\|_{L^1(K)}~ \lambda^{-1}.
\end{align}

Now we just have to handle the term 
\begin{align*}
\bigg|\bigg(\bigcup_iW_i\bigg)^\complement \bigcap \{x \in K : |Bf_1(x)|> \lambda/2\} \bigg|,
\end{align*}
\ni for this one we prove that 
\begin{align} \label{eq11}
\text{if}~  x \notin \cup_iW_i ~ \text{then}~ Bf_1(x) =0,
\end{align}
\ni and which gives that the above term is zero and reduces (\ref{eq9}) to one term only. We turn to proving (\ref{eq11}), 
\begin{align*}
 g_j \ast f_1(x) & = \int\limits_{K} a(y)\Phi_{A_j}(y)f_1(x-y)dy\\
&= \int\limits_{\mathfrak{D}^*} a(y)f_1(x-y)dy,
\end{align*}

\ni since $\text{supp}~ a\Phi_{A_j} \subset \mathfrak{D}^* \cap A_j = \mathfrak{D}^*.$ Thus 

\begin{align*}
g_j \ast f_1(x) &= \int\limits_{y \in (x - W_i) \cap \mathfrak{D}^*} a(y)f_1(x-y)dy.
\end{align*}

Let us note that $0 \notin x - W_i.$ If  $0 \in x - W_i$ then  $x \in W_i,$  a contradiction since $x \in (\cup_iW_i)^\complement.$ Then $(x - W_i) \cap \mathfrak{D}^* = \emptyset$ and $g_j \ast f_1(x)$ is zero. Hence 

\begin{align*} 
Bf_1(x) & = \sum\limits_{j = k }^{\infty }q^{-(j+1)} g_j \ast f_1(x)\\
& = 0.
\end{align*}

By (\ref{eq7}), (\ref{eq8}) and (\ref{eq10}), we obtain

\begin{align} \nonumber
|\{x \in K : |Bf(x)| >  \lambda\}| & \leq \lambda^{-1}~ \|f\|_{L^1(K)} + 4q \lambda^{-1}~ \|f\|_{L^1(K)} \\ \label{eq12}
& = (1+4q)\lambda^{-1}~\|f\|_{L^1(K)} .
\end{align}

This concludes the proof of the fact that $B$ is weak type $(1,1).$ From (\ref{eq6}) and (\ref{eq12}) we have

\begin{align*}
\|Bf\|_{L^2(K)} \leq q^{-k}(q-1)^{-1} \|f\|_{L^2(K)},
\end{align*}
and 
  \begin{align*}
  |\{x \in K : |Bf(x)| >  \lambda\}| \leq (1+4q) \lambda^{-1} \|f\|_{L^1(K)}.
  \end{align*}

Applying Lemma \ref{lemma2} we obtain that $T_k$ is of type $(r,r)(1 < r \leq 2),$ i.e, for $f \in S(K),$ $Bf \in L^r(K)(1 < r \leq 2)$ and

\begin{align} \label{n1}
\|Bf\|_{L^r(K)} \leq Aq^{-k} \|f\|_{L^r(K)},
\end{align}

where $A$ is independent of $k.$


\subsection{$B$ commutes with translations} Let $h \in K$ and $f \in S(K),$
\begin{align*}
\widehat{\tau_h(Bf)} & = \sum\limits_{j = k }^{\infty }q^{-(j+1)} \widehat{(\tau_h (g_j \ast f))}\\ 
& = \sum\limits_{j = k }^{\infty }q^{-(j+1)} \bar{\chi_h} \widehat{(g_j \ast f)}~~~~~~~~~(\because~ \widehat{\tau_hf}=\bar{\chi_h}\widehat{f}) \\ 
& = \sum\limits_{j = k }^{\infty }q^{-(j+1)} \widehat{g_j} (\bar{\chi_h} \widehat{ f}) \\
& = \sum\limits_{j = k }^{\infty }q^{-(j+1)} \widehat{g_j} ( \widehat{ \tau_h f}) \\
& = \sum\limits_{j = k }^{\infty }q^{-(j+1)} \widehat{(g_j \ast \tau_hf)} \\
& = \widehat{(B\tau_hf)}, 
\end{align*}

 which implies that $B$ commutes with translations. From (\ref{n1}) and Lemma \ref{lemma3}, we obtain that for $f \in S(K),$ $Bf \in L^r(K)(2 < r < \infty)$ and

\begin{align} \label{eq14}
\|Bf\|_{L^r(K)} \leq  Aq^{-k} \|f\|_{L^r(K)}. 
\end{align}

Combining (\ref{n1}) and (\ref{eq14}), we get whenever $f \in S(K),$ $Bf \in L^r(K)$ and

\begin{align*} 
\|Bf\|_{L^r(K)} \leq  Aq^{-k} \|f\|_{L^r(K)},~~~~~ 1 < r < \infty. \\
\end{align*}

Therefore by (\ref{eq3}), for  $f \in S(K),$ $T_kf \in L^r(K)$ and

\begin{align} \label{neq1}
\|T_kf\|_{L^r(K)} \leq  Cq^{-k} \|\Omega\|_{H^1(\mathfrak{D}^*)}  \|f\|_{L^r(K)},~~~~~ 1 < r < \infty. 
\end{align}

Let $f \in L^r(K),$ since $S(K)$ is dense in $L^r(K)~1 \leq r < \infty,$ i.e. there is a sequence $(f_j)_{j \in \mathbb{N}} \subseteq S(K)$ with $f_j \rightarrow f$ in $L^r(K).$ From (\ref{neq1}), we see that 

\begin{align*}
\|T_kf_j - T_kf_{j^\prime}\|_{L^r(K)} & = \|T_k(f_j - f_{j^\prime})\|_{L^r(K)} \\
& \leq   Cq^{-k} \|\Omega\|_{H^1(\mathfrak{D}^*)}  \|f_j - f_{j^\prime}\|_{L^r(K)},  
\end{align*}

and so $(T_kf_j)_{j \in \mathbb{N}} $ is a Cauchy sequence in $L^r(K).$ By completeness, the limit  $T_kf_j \rightarrow T_kf$ exists in $L^r(K)$ and by (\ref{neq1}),
\begin{align*}
\|T_kf\|_{L^r(K)} &= \lim_{j\to\infty}\|T_kf_j\|_{L^r(K)} \\
& \leq  Cq^{-k} \|\Omega\|_{H^1(\mathfrak{D}^*)} \lim_{j\to\infty} \|f_j\|_{L^r(K)} \\
 & \leq  Cq^{-k} \|\Omega\|_{H^1(\mathfrak{D}^*)} \|f\|_{L^r(K)},
\end{align*}
which completes the proof of Theorem \ref{thm1}.



\section{Proof of Theorem \ref{thm2}}

Let $1 < r, t < \infty$ and $s > 0.$ Let $\Omega \in H^1(\mathfrak{D}^*)$ then it follows from (\ref{eq3}) that
\begin{align*}
T_k(f)(x) & = \sum\limits_{i=1}^{\infty} \lambda_i Bf, 
\end{align*}
and
\begin{align} \nonumber
\|T_kf\|_{F_{r t}^s(K)} & \leq  \sum\limits_{i=1}^{\infty} |\lambda_i| \|Bf\|_{F_{r t}^s(K)} \\ \label{eq15} 
& \leq C \|\Omega\|_{H^1(\mathfrak{D}^*)} \|Bf\|_{F_{r t}^s(K)}, \\ \nonumber
\end{align}
\ni where
\begin{align*}
Bf(x) = \sum\limits_{j = k }^{\infty }q^{-(j+1)} g_j \ast f(x).
\end{align*}
Then 
\begin{align*}
\|Bf\|_{F_{r t}^s(K)} \leq  \sum_{j = k }^{\infty } q^{-(j+1)}~  \|g_j \ast f(x) \|_{F_{r t}^s(K)}.
\end{align*}

So, to estimate the $F_{r t}^s(K)$ norm of $T_k$ it is enough to find $\|g_j \ast f(x) \|_{F_{r t}^s(K)}.$ We begin with the observations:
\begin{align} \nonumber
\Delta_l(g_j \ast f) & = \Delta_l(f \ast g_j) \\ \nonumber
& =  \mathcal{F}^{-1} \phi_l \mathcal{F}(f \ast g_j) \\
& =  (\mathcal{F}^{-1} \phi_l \mathcal{F}f) \ast  g_j, 
\end{align}

where $\phi_l$ is decomposition of unity see Definition \ref{def2} and 

\begin{align}
\|(\mathcal{F}^{-1} \phi_l \mathcal{F}f) \ast g_j\|_{\ell^t} \leq (\|\mathcal{F}^{-1} \phi_l \mathcal{F}f \|_{\ell^t}) \ast g_j,
\end{align}

which follows from (\ref{8}) with $F_l= (\mathcal{F}^{-1} \phi_l \mathcal{F}f) \ast g_j.$ Now

\begingroup
\addtolength{\jot}{.5em}
\begin{align*}
\|g_j \ast f(x) \|_{F_{r t}^s(K)} & = \Bigg\|\Bigg\{\sum_{l=0}^{\infty} |q^{sl} \Delta_l(g_j \ast f)|^t   \Bigg\}^{\frac{1}{t}} \Bigg\|_{L^r(K)} \\
& = \Bigg\|\Bigg\{\sum_{l=0}^{\infty} q^{slt} |(\mathcal{F}^{-1} \phi_l \mathcal{F}f) \ast g_j|^t   \Bigg\}^{\frac{1}{t}} \Bigg\|_{L^r(K)} \\
& \leq  \Bigg\|\Bigg\{\sum_{l=0}^{\infty} q^{slt} | \mathcal{F}^{-1} \phi_l  \mathcal{F}f |^t   \Bigg\}^{\frac{1}{t}} \ast g_j \Bigg\|_{L^r(K)} \\
& \leq \Bigg\|\Bigg\{\sum_{l=0}^{\infty} q^{slt} | \mathcal{F}^{-1} \phi_l  \mathcal{F}f |^t   \Bigg\}^{\frac{1}{t}}\Bigg\|_{L^r(K)} \|g_j\|_{L^1(K)} \\
& \leq  \|f\|_{F_{r t}^s(K)}. 
\end{align*}
\endgroup

Therefore
\begin{align} \nonumber
\|Bf\|_{F_{r t}^s(K)} & \leq  \sum_{j = k }^{\infty } q^{-(j+1)}~ \|f\|_{F_{r t}^s(K)} \\ \label{eq16}
& =  q^{-k}(q-1)^{-1} \|f\|_{F_{r t}^s(K)}. \\ \nonumber
\end{align}

By (\ref{eq15}) and (\ref{eq16}), we have

\begin{align*}
\|T_kf\|_{F_{r t}^s(K)} & \leq C  \|\Omega\|_{H^1(\mathfrak{D}^*)} \|Bf\|_{F_{r t}^s(K)} \\
 & \leq C q^{-k} \|\Omega\|_{H^1(\mathfrak{D}^*)} \|f\|_{F_{r t}^s(K)}.\\
\end{align*}

Since $S(K)$ is dense in $F_{r t}^s(K)~1 \leq r,t < \infty,$ therefore for $f \in F_{r t}^s(K)$ we have

\begin{align*}
\|T_kf\|_{F_{r t}^s(K)} & \leq Cq^{-k} \|\Omega\|_{H^1(\mathfrak{D}^*)} \|f\|_{F_{r t}^s(K)}.\\
\end{align*}

This is the conclusion ($\romannumeral 1$) of Theorem \ref{thm2}. By Theorem \ref{thm1}, we have for $1 < r,t < \infty$ 

\begin{align*}
\bigg(\sum_{l=0}^{\infty} \|T_kf_j\|^t_{L^r(K)} \bigg)^{1/t} \leq  Cq^{-k} \|\Omega\|_{H^1(\mathfrak{D}^*)} \bigg(\sum_{l=0}^{\infty} \|f_j\|^t_{L^r(K)} \bigg)^{1/t}. \\
\end{align*}

Then we get for any $1 < r,t < \infty$  and $s > 0,$

\begin{align*}
 \|T_kf\|_{B_{r t}^s(K)} &= \Bigg(\sum_{l=0}^{\infty} q^{slt}\| \Delta_l(T_kf) \|_{L^r(K)}^t   \Bigg)^{\frac{1}{t}} \\ 
& = \Bigg(\sum_{l=0}^{\infty} q^{slt} \| T_k(\Delta_lf) \|_{L^r(K)}^t   \Bigg)^{\frac{1}{t}} \\ 
& \leq   Cq^{-k}  \|\Omega\|_{H^1(\mathfrak{D}^*)} \Bigg(\sum_{l=0}^{\infty}  q^{slt} \|\Delta_lf \|_{L^r(K)}^t   \Bigg)^{\frac{1}{t}} \\ 
 & =   Cq^{-k}  \|\Omega\|_{H^1(\mathfrak{D}^*)} \|f\|_{B_{r t}^s(K)}. \\ 
\end{align*}

Again using the fact that $S(K)$ is dense in $B_{r t}^s(K)~1 \leq r,t < \infty,$ therefore for $f \in B_{r t}^s(K)$ we have 

\begin{align*}
\|T_kf\|_{B_{r t}^s(K)} & \leq  Cq^{-k}  \|\Omega\|_{H^1(\mathfrak{D}^*)} \|f\|_{B_{r t}^s(K)}.\\
\end{align*}

This concludes the proof of Theorem \ref{thm2}.

\bibliographystyle{achemso}

\end{document}